\documentclass[12pt]{article}
\usepackage{tikz}
\usepackage{scalefnt}
\usepackage{multicol}
\usepackage{enumerate}

\usepackage[latin1]{inputenc}
\usepackage[T1]{fontenc}
\usepackage{amsmath,amssymb,amsthm}
\usepackage{graphicx}
\usepackage{mathptmx}
\usepackage{enumitem}

\newtheorem{theorem}{Theorem}[section]

\theoremstyle{definition}

\newcommand{\R}{\mathbb{R}}
\newcommand{\Q}{\mathbb{Q}}

\usepackage[total={5.5in,8.5in},text={4.7in,7.65in},centering]{geometry}

\begin{document}
\pagenumbering{gobble} 



{\bf \Large \centering  NEW LOWER BOUNDS ON $\chi(\mathbb{R}^d)$ FOR $d=8, \dots, 12$\\}

\bigskip

{ \centering \Large \bf Matthew Kahle\\}
{ \centering \large Ohio State University\\}
{ \centering \large mkahle@math.osu.edu\\}
{ \centering \Large \bf Birra Taha\\}
{\centering \large Cornell University\\}
{ \centering \large bit2002@med.cornell.edu\\}

\bigskip

%
%

\section{Introduction}

We study the infinite graph whose vertices correspond to points in $d$-dimensional Euclidean space and where two vertices $p,q$ are considered adjacent $p \sim q$ if $d(p,q) = 1$. We denote this graph by $\R^d$.

Recall that the chromatic number $\chi(G)$ of a graph $G$ is the smallest number of colors needed for a proper coloring, i.e.\ the smallest $k$ so that there exists a function $f:V(G) \to \{1, 2, \dots, k \}$ such that whenever $p \sim q$, we have $f(p) \neq f(q)$.

Determining the chromatic number of the plane $\chi \left( \R^2 \right)$ is a well known problem in combinatorics. For several decades the bounds have held fast at 
$$ 4 \le \chi \left( \R^2 \right) \le 7.$$ See Soifer's ``The mathematical coloring book'' \cite{Soifer} for an encyclopedic overview of this problem and its history.

The problem has also been studied in higher dimensions. See \cite{Soifer} and Kupavskii--Raigorodskii's paper \cite{KR09}. With the aid of extensive calculations in the free and open-source software Sage \cite{sage}, we obtain new lower bounds for $\chi( \R^d)$ in dimensions $d = 8 \dots 12$. We verified all of the computations in Maple 17.

\newpage

Our main result is the following.

\begin{theorem}[]
We have that
$\chi \left( \R^8 \right) \ge 19$, 
$\chi \left( \R^9 \right)  \ge 21$,
$\chi \left( \R^{10} \right)  \ge 26$,
$\chi \left( \R^{11} \right)  \ge 32$, and
$\chi \left( \R^{12} \right) \ge 32.$
\end{theorem}

As far as we are aware, the best previously published lower bounds were, respectively:

\noindent $\chi \left( \R^8 \right) \ge 16$ by Larman and Rogers in 1972 \cite{LR72}, and  
$\chi \left( \R^9 \right)  \ge 21$,
$\chi \left( \R^{10} \right)  \ge 23$,
$\chi \left( \R^{11} \right)  \ge 23$, and
$\chi \left( \R^{12} \right) \ge 25$ by Kupavskii and Raigorodskii in 2009 \cite{KR09}.


\medskip

The rest of the paper is organized as follows. We prove the bounds $\chi \left( \R^{10} \right)  \ge 26$, $\chi \left( \R^{11} \right)  \ge 32$, and
$\chi \left( \R^{12} \right) \ge 32$ in Section \ref{sec:hcube}; the bound $\chi \left( \R^9 \right)  \ge 21$ in Section \ref{sec:hsec}; and the bound $\chi \left( \R^8 \right) \ge 19$ in Section \ref{sec:E8}.

\medskip

{\bf Acknowledgements.} We thank the Ohio Supercomputer Center for computing resources and technical support. This work was sponsored in part by the first author's Sloan Research Fellowship. We thank Chris Peterson for helpful and inspiring conversations, and in particular for pointing out that the $8$-dimensional Gosset polytope already gives the lower bound 
$\chi \left( \R^8 \right) \ge 15$.

\section{Hypercube graphs} \label{sec:hcube}

Graphs constructed on vertices of the $d$-dimensional cube $ \{ 0, 1 \}^d$ provide important examples in geometric graph theory. Frankl and Wilson's proof that
$$\chi \left( \R^d \right) \ge \exp(C d)$$
for some constant $C>0$, for example, uses such graphs \cite{FW81}. See also Kahn and Kalai's subsequent counterexample to Borsuk's conjecture \cite{KK93}.

Define the hypercube graph $C(d,u)$ to have vertices $V= \{ 0, 1 \}^d$, with two vertices adjacent if their Hamming distance is $u$. Note that Hamming distance $u$ corresponds to Euclidean distance $\sqrt{u}$. So by dilating Euclidean space by a factor of $1 / \sqrt{u}$, we see that $C(d,u)$ is a unit-distance graph in $\R^d$.

For $d$ odd $C(d,u)$ is bipartite, so for our purposes these are not very interesting. For $d$ even $C(d,u)$ has two isomorphic connected components --- we denote one of these ``half-cube'' connected components on $2^{d-1}$ vertices  by $H(d,u)$. In general we have that
$$\chi(R^d) \ge \chi(C(d,u)) =  \chi (H(d,u)).$$

For example $C(5,2)$ is a graph on $32$ vertices, regular of degree $10$, and the half-cube $H(5,2)$ is a graph on $16$ vertices.

Recall that for a graph $G$ with $|V(G)|$ vertices and independence number $\alpha(G)$, we have that
$$\chi(G) \ge \frac{|V(G)|}{\alpha(G)}.$$

One checks that the independence number $\alpha(H(5,2)) = 2$. Then the independence-number bound gives that $\chi (H) \ge 8$,
and since $H$ is a unit-distance graph in $\R^5$, this also gives that $\chi(\R^5) \ge 8$. This half-cube example is well known; see for example \cite{LR72}.\\

\medskip

In Figure \ref{fig:asy} we summarize our results on chromatic numbers of hypercube graphs for small $d$ and $u$. In some cases we were able to compute $\chi( C(d,u) )$ exactly. In some other cases we were not, but we were still able to compute the independence number $\alpha ( C(d,u) )$, giving a lower bound.

\begin{figure}
\scalefont{2}
\scalebox{0.99}{
\begin{tikzpicture}[thick]

\hspace{-.1in}

\large


\draw (1,1) -- (11,1);
\draw (1,2) -- (11,2);
\draw (1,3) -- (11,3);
\draw(1,4) -- (11,4);

\draw (1,1) -- (1,4);
\draw (2,1) -- (2,4);
\draw (3,1) -- (3,4);
\draw (4,1) -- (4,4);
\draw (5,1) -- (5,4);
\draw (6,1) -- (6,4);
\draw (7,1) -- (7,4);
\draw (8,1) -- (8,4);
\draw (9,1) -- (9,4);
\draw (10,1) -- (10,4);
\draw (11,1) -- (11,4);



\node at (6,-1/2) {\bf  $d$};
\node at (-1/2,5/2) {\bf  $u$};

\node at (3/2,1/2) {\bf 2};
\node at (5/2,1/2) {\bf 3};
\node at (7/2,1/2) {\bf 4};
\node at (9/2,1/2) {\bf 5};
\node at (11/2,1/2) {\bf 6};
\node at (13/2,1/2) {\bf 7};
\node at (15/2,1/2) {\bf 8};
\node at (17/2,1/2) {\bf 9};
\node at (19/2,1/2) {\bf 10};
\node at (21/2,1/2) {\bf 11};

\node at (1/2,3/2) {\bf 2};
\node at (1/2,5/2) {\bf 4};
\node at (1/2,7/2) {\bf 6};

\small

%

\node at (3/2, 3/2) {$2$};
\node at (5/2, 3/2) {$4$};
\node at (7/2, 3/2) {$4$};
\node at (9/2, 3/2) {$8$};
\node at (11/2, 3/2) {$8$};
\node at (13/2, 3/2) {$8$};
\node at (15/2, 3/2) {$8$};
\node at (17/2, 3/2) {$\ge{13}$};
\node at (19/2, 3/2) {$\ge{13}$};
\node at (21/2, 3/2) {$\ge{13}$};

%

\node at (3/2, 5/2) {$1$};
\node at (5/2, 5/2) {$1$};
\node at (7/2, 5/2) {$2$};
\node at (9/2, 5/2) {$4$};
\node at (11/2, 5/2) {$7$};
\node at (13/2, 5/2) {$8$};
\node at (15/2, 5/2) {$8$};
\node at (17/2, 5/2) {$\ge{15}$};
\node at (19/2, 5/2) {$\ge{26}$};
\node at (21/2, 5/2) {$\ge{32}$};
%
%

\node at (3/2, 7/2) {$1$};
\node at (5/2, 7/2) {$1$};
\node at (7/2, 7/2) {$1$};
\node at (9/2, 7/2) {$1$};
\node at (11/2, 7/2) {$2$};
\node at (13/2, 7/2) {$4$};
\node at (15/2, 7/2) {$\ge{5}$};
\node at (17/2, 7/2) {$\ge{7}$};
\node at (19/2, 7/2) {$\ge{11}$};
\node at (21/2, 7/2) {$\ge{16}$};

%

%
%
%
%

\end{tikzpicture}
}
\caption{Chromatic numbers of some hypercube graphs $\chi (C(d,u))$ for some small $d$ and $u$. We restrict to $u$ even, since $C(d,u)$ is bipartite for $u$ odd. In some cases we compute chromatic number exactly, and have an independence-number lower bound. For $(d,u) = (10,2) \mbox{ and } (11,2)$ the lower bound is by monotonicity of rows.}
\label{fig:asy}
\end{figure}

We obtain new lower bounds on $\chi( \R^d )$ for $d = 10$ and $d=11$ this way, in both cases with Hamming distance $u=4$. According to a calculation in Sage, the halfcube $H(10,4)$ has independence number $\alpha( H(10,4) ) = 20$. There are $512$ vertices, so then 
$$ \chi \left( H(10,4) \right) \ge \frac{512}{20} = 25.6,$$
and since the chromatic number is an integer,
$$ \chi \left( \R^{10} \right) \ge 26.$$

According to another calculation in Sage, the half-cube $H(11,4)$ has independence number $\alpha( H(10,4) ) = 32$. There are $1024$ vertices, and then we have 
$$ \chi \left( R^{11} \right) \ge \frac{1024}{ 32} = 32.$$
This last example also gives a new record lower bound $\chi( \R^{12}) \ge 32$.

Since Hamming distance $4$ corresponds to Euclidean distance $2$, these graphs still have rational coordinates when rescaled to unit-distance graphs. 

\section{Hyperplane slices of hypercube graphs} \label{sec:hsec}

Let $C(d,u,s)$ denote the intersection of $C(d,u)$ with the hyperplane at height $s$
$$x_1 + \dots +  x_d = s.$$
Clearly, $C(d,u,s)$ is a unit distance graph in $R^{d-1}$, so 
$$\chi \left( R^{d-1} \right) \ge \chi \left( C(d,u,s) \right).$$ 

For example, $C(10,4)$ has $2^{10}=1024$ vertices and is regular of degree ${10 \choose 4}$, where the subgraph $C(10,4,5)$ has ${10 \choose 5} = 252$ vertices and is regular of degree $100$.

We found with a computation in Sage that $$\alpha \left( C(10,4,5) \right) = 12,$$
so $$\chi \left( R^9 \right) \ge \frac{252}{12} = 21.$$

This gives an alternate proof for the currently best known bound of Kupavskii and Raigorodskii \cite{KR09}.

\section{A unit-distance graph in $\R^8$} \label{sec:E8}
Let $G_0$ be a graph whose vertices are the $240$ shortest vectors in the $E8$ lattice. 
These may also be described as the vertices of an $8$-dimensional Gosset polytope.

Concretely, there are $112$ vertices with integer entries obtained from arbitrary permutations of the vectors
$$ ( \pm 2, \pm 2, 0,0,0,0,0,0),$$
and $128$ vertices with integer entries obtained from all vectors
$$ (\pm 1, \pm 1, \pm 1, \pm 1, \pm 1, \pm 1, \pm 1, \pm 1) $$ with an even number of minus signs.

Adjacency is with respect to Euclidean distance $4$. So for example,
$ (2, 2, 0,0,0,0,0,0)$ is adjacent to $(0,0, 2, 2, 0,0,0,0)$.

The graph $G_0$ lives on a $7$-dimensional sphere of radius $2 \sqrt{2}$. A pair of vertices on this sphere is adjacent if and only if the corresponding vectors are orthogonal.

Now define $P$ to be the set of integer points within distance $4$ of the origin.
%

\medskip

Our algorithm is as follows.

\begin{enumerate}[noitemsep]
\item Initialize $G_0$ as above.
\item Choose a vertex $x \in P$ outside of the current graph $G_i$.
\item If the independence number remains unchanged on adding vertex $x$ to the graph $G_i$ then add $X$, i.e.\ if $\alpha(G_i + x) = \alpha(G_i)$, then set $V(G_{i+1} ):= V(G_i) + x$.
\item If there are any more points in $P$ outside of the graph $V(G_i)$ , go to step (2).
\end{enumerate}

After experimentation, we found that at least $49$ points (Figure \ref{fig:add}) could be added without increasing the independence number, which then gives a unit-distance graph $G$ with 
$$\chi(G) \ge \frac{289}{16} = 18.0625,$$
so $\chi(\R^8) \ge 19$.  Our graph $G$ can be rescaled to have rational coordinates and unit distance, so this also gives the new lower bound $\chi(\Q^8) \ge 19$.

\bigskip

\begin{figure}
\begin{multicols}{2}
\begin{itemize}[noitemsep]
\item (-2, 0, -2, 0, 0, 2, 0, 2)
\item(-2, 0, 0, 0, 2, 0, -2, 2)
\item (-2, 0, 0, 2, 2, 2, 0, 0)
\item (-2, 0, 2, 0, 0, 2, 0, 2)
\item (-2, 0, 2, 2, 0, 0, 0, 2)
\item (-2, 2, 0, 0, -2, 2, 0, 0)
\item (-2, 2, 0, 0, 0, 0, 2, 2)
\item (-2, 2, 2, 0, 2, 0, 0, 0)
\item (0, -2, -2, 0, 0, 2, 0, 2)
\item (0, -2, -2, 0, 2, 0, 0, 2)
\item (0, -2, 0, 2, 2, 0, 0, 2)
\item (0, 0, -2, 0, 2, 0, -2, 2)
\item (0, 0, -2, 0, 2, 2, 2, 0)
\item (0, 0, -2, 2, 2, 0, 2, 0)
\item (0, 0, -1, 1, -1, 0, 0, 1)
\item (0, 0, 0, 2, 0, 2, 2, 2)
\item (0, 0, 0, 2, 2, 0, -2, 2)
\item (0, 0, 2, 0, -2, 0, 2, 2)
\item (0, 0, 2, 0, 2, -2, -2, 0)
\item (0, 0, 2, 2, 2, 0, 2, 0)
\item (0, 2, 0, 0, -2, -2, 0, 2)
\item (0, 2, 0, 0, 2, 0, -2, -2)
\item (0, 2, 0, 2, 0, 2, 0, -2)
\item (0, 2, 0, 2, 2, 2, 0, 0)
\item (1, -1, 1, -1, 1, 1, -1, 3)
\item (1, 0, 0, 1, 1, 1, 0, 0)
\item (1, 0, 1, 0, 1, 1, 0, 0)
\item (1, 1, 0, 1, 1, 0, 0, 0)
\item (1, 1, 1, 0, 0, 0, 0, 1)
\item (1, 1, 1, 1, 0, 0, 0, 0)
\item (2, -2, -2, 0, 0, 0, 2, 0)
\item (2, -2, -2, 0, 0, 2, 0, 0)
\item (2, -2, 0, -2, 0, 0, 2, 0)
\item (2, -2, 0, 0, -2, 0, 0, 2)
\item (2, -2, 0, 0, 2, 0, -2, 0)
\item (2, 0, 0, 0, 0, -2, 2, -2)
\item (2, 0, 0, 0, 0, 2, 2, -2)
\item (2, 0, 0, 0, 2, -2, 0, -2)
\item (2, 0, 0, 0, 2, 2, 0, -2)
\item (2, 0, 0, 0, 2, 2, 2, 0)
\item (2, 0, 0, 2, 0, -2, 2, 0)
\item (2, 0, 0, 2, 2, 0, 0, -2)
\item (2, 0, 0, 2, 2, 0, 0, 2)
\item (2, 0, 2, -2, -2, 0, 0, 0)
\item (2, 0, 2, 0, 0, 0, 2, -2)
\item (2, 0, 2, 2, 0, 0, 0, 2)
\item (2, 2, 0, 0, 0, 0, 2, 2)
\item (3, -1, 1, -1, 1, -1, 1, 1)
\item (3, -1, 1, 1, -1, -1, -1, -1)
\end{itemize}
\end{multicols}
\caption{List of $49$ additional points added to the vertices of the $8$-dimensional Gosset polytope to obtain a unit-distance graph with chromatic number at least $19$.}
\label{fig:add}
\end{figure}

Mann used extensive computer-aided calculations to hunt for unit-distance graphs with large chromatic number in 2003 \cite{M03} --- he established the lower bounds $\chi(\Q^6) \ge 10$, $\chi(\Q^7) \ge 13$, and $\chi(\Q^8) \ge 16$.

Cibulka studied similar constructions in \cite{Cibulka08}. He also added points to Gosset polytopes one at a time, and used computer-aided calculations to establish the bounds $\chi(\Q^5) \ge 8$ and $\chi(\Q^7) \ge 15$.

Unit-distance graphs built on Gosset polytopes go back at least to the landmark paper of Larman and Rogers in 1972 \cite{LR72}. They in turn thanked McMullen for suggesting the idea. Larman and Rogers exhibited a configuration of $64$ points, a ``spindle'' over the $7$-dimensional Gosset polytope, to establish the previous record lower bound $$\chi( \R^8) \ge 16.$$

\bigskip

\bibliographystyle{plain}
\bibliography{lbcrefs}

\begin{thebibliography}{1}

\bibitem{Cibulka08}
Josef Cibulka.
\newblock On the chromatic number of real and rational spaces.
\newblock {\em Geombinatorics}, 18(2):53--65, 2008.

\bibitem{FW81}
P.~Frankl and R.~M. Wilson.
\newblock Intersection theorems with geometric consequences.
\newblock {\em Combinatorica}, 1(4):357--368, 1981.

\bibitem{KK93}
Jeff Kahn and Gil Kalai.
\newblock A counterexample to {B}orsuk's conjecture.
\newblock {\em Bull. Amer. Math. Soc. (N.S.)}, 29(1):60--62, 1993.

\bibitem{KR09}
A.~B. Kupavskii and A.~M. Raigorodskii.
\newblock On the chromatic numbers of small-dimensional {E}uclidean spaces.
\newblock In {\em European {C}onference on {C}ombinatorics, {G}raph {T}heory
  and {A}pplications ({E}uro{C}omb 2009)}, volume~34 of {\em Electron. Notes
  Discrete Math.}, pages 435--439. Elsevier Sci. B. V., Amsterdam, 2009.

\bibitem{LR72}
D.~G. Larman and C.~A. Rogers.
\newblock The realization of distances within sets in {E}uclidean space.
\newblock {\em Mathematika}, 19:1--24, 1972.

\bibitem{M03}
Matthias Mann.
\newblock Hunting unit-distance graphs in rational {$n$}-spaces.
\newblock {\em Geombinatorics}, 13(2):86--97, 2003.

\bibitem{Soifer}
Alexander Soifer.
\newblock {\em The mathematical coloring book}.
\newblock Springer, New York, 2009.
\newblock Mathematics of coloring and the colorful life of its creators, With
  forewords by Branko Gr{\"u}nbaum, Peter D. Johnson, Jr. and Cecil Rousseau.

\bibitem{sage}
W.\thinspace{}A. Stein et~al.
\newblock {\em {S}age {M}athematics {S}oftware ({V}ersion 5.10)}.
\newblock The Sage Development Team, 2013-06-17.
\newblock {\tt http://www.sagemath.org}.

\end{thebibliography}

\end{document}